\input ./style/arxiv-general.cfg
\documentclass[aos,MSNbibl,seceqn,dvips]{arximspdf}
\makeatletter
   \@ifpackageloaded{graphicx}{}{\usepackage{graphicx}}
\makeatother

\doi{10.1214/15-AOS1412}
\volume{44}
\issue{4}
\pubyear{2016}
\firstpage{1438}
\lastpage{1466}
\docsubty{FLA}

\makeatletter
\newcommand{\lleft}{\left}
\newcommand{\rrvert}{\vert}
\newcommand{\rright}{\right}
\newcommand{\llvert}{\vert}
\def\sfrac#1#2{#1/#2}
\def\vfrac#1#2{(#1)/#2}
\def\afrac#1#2{#1/(#2)}

\newcommand{\eqref}[1]{(\ref{#1})}
\let\epsilon\varepsilon
\newtheorem{theorem}{Theorem}[section]
\newtheorem{lemma}[theorem]{Lemma}

\newproclaim{definition}[theorem]{Definition}
\newproclaim{conjecture}[theorem]{Conjecture}
\newproclaim{remark}{Remark}[section]
\newproclaim{remarks}{Remarks}[section]
\newproclaim{example}{Example}[section]
\makeatother

\begin{document}
\begin{frontmatter}

\title{Global testing against sparse alternatives in~time-frequency analysis\thanksref{T1}}
\runtitle{Signal Detection}

\begin{aug}
\author[A]{\fnms{T. Tony}~\snm{Cai}\ead[label=e1]{tcai@wharton.upenn.edu}},
\author[B]{\fnms{Yonina C.}~\snm{Eldar}\ead[label=e2]{yonina@ee.technion.ac.il}}
\and
\author[C]{\fnms{Xiaodong}~\snm{Li}\corref{}\ead[label=e3]{xdgli@ucdavis.edu}}
\runauthor{T. T. Cai, Y. C.  Eldar and X. Li}
\affiliation{University of Pennsylvania,  Technion---Israel Institute of Technology\\ and University of California, Davis}
\address[A]{T. T. Cai\\
Department of Statistics\\
The Wharton School\\
University of Pennsylvania\\
400 Jon M. Huntsman Hall\\
3730 Walnut Street\\
Philadelphia, Pennsylvania 19104-6340\\
USA\\
\printead{e1}}
\address[B]{Y. C.  Eldar\\
Department of Electrical Engineering\\
Technion---Israel Institute\\
\quad of Technology\\
Haifa 32000\\
Israel\\
\printead{e2}}
\address[C]{X. Li\\
Department of Statistics\\
University of California, Davis\\
4109 Mathematical Sciences\\
Davis, California 95616\\
USA\\
\printead{e3}}
\end{aug}
\thankstext{T1}{Supported by  NSF Grants
DMS-12-08982 and DMS-14-03708, NIH Grant R01 CA127334, and the Wharton
Dean's Fund for Post-Doctoral Research.}

\received{\smonth{6} \syear{2015}}
%
\revised{\smonth{10} \syear{2015}}

%
\begin{abstract}
In this paper, an \emph{over-sampled periodogram higher criticism}
(OPHC) test is proposed for the global detection of sparse periodic
effects in a complex-valued time series. An explicit minimax detection
boundary is established between the rareness and weakness of the
complex sinusoids hidden in the series. The OPHC test is shown to be
asymptotically powerful in the detectable region. Numerical simulations
illustrate and verify the effectiveness of the proposed test.
Furthermore, the periodogram over-sampled by $O(\log N)$ is proven
universally optimal in global testing for periodicities under a mild
minimum separation condition.
\end{abstract}

%
\begin{keyword}[class=AMS]
\kwd{62C20}
\kwd{62F03}
\kwd{62F05}
\kwd{62G30}
\kwd{62G32}
\end{keyword}
\begin{keyword} 
\kwd{Testing for periodicity}
\kwd{sparsity}
\kwd{over-sampled periodogram}
\kwd{higher criticism}
\kwd{detection boundary}
\kwd{empirical processes}
\end{keyword}
\end{frontmatter}

\section{Introduction}
In this paper, we study the problem of global testing for periodicity.
Suppose $u_t, t=1, \ldots, N$, is a real-valued time series observed
at equispaced time points, that satisfies the model
%
\begin{equation}
\label{eqrealmodel} u_t=\sum_{j=1}^s
a_j \sin(\omega_j t + \phi_j) +
\epsilon_t,
\end{equation}
where the noise $\epsilon_t \sim\mathcal{N}(0, \sigma^2)$ are
i.i.d. normal variables. In the complex-valued case, similarly, the
observed series satisfies the model
%
\begin{equation}
\label{eqcomplexmodel} y_t=\sum_{j=1}^s
a_j e^{i(\omega_j t + \phi_j)} + z_t,
\end{equation}
where $z_t$ represents zero-mean i.i.d. complex white noise, that is,
$z_t=z_{1t} + iz_{2t}$ with $
[ z_{1t}, z_{2t} ]
\sim\mathcal{N}(\mathbf{0}, \frac{\sigma^2}{2}\mathbf{I}_{2})$. In both
cases, we are interested in testing
%
\begin{equation}
\label{eqtesting} H_0: a_1=\cdots=a_s=0 \quad
\mbox{versus}\quad H_1: a_j >0\qquad  \mbox{for at least one
$j$},
\end{equation}
in which periodicity exists in the series under the alternative.

Global detection of periodic patterns in time series analysis have
various applications. We give several examples as follows:


\subsection*{Signal detection} Global detection of sinusoidal signals
is a fundamental signal processing task prior to information extraction
\cite{VanTrees92,Whalen95,Kay98}. As summarized in \cite{Kay98},
global testing for waveforms can be utilized in radar and sonar
systems, such as detecting whether aircrafts are approaching \cite
{Skolnik80} and detecting whether enemy submarines are present \cite{KPK1981}.

\subsection*{Gene expression studies} Global detection of periodic
patterns in time series due to various biological rhythms such as cell
division, circadian rhythms, life cycles of microorganisms and many
others is an important problem in gene expression studies \cite
{WFS04,Chen05,ALPHY05,GCM06,ALY07}. For example, in order to
identify a collection of genes which are responsible in the cell cycle,
it was proposed in \cite{WFS04} to implement global testing of
periodicity for each gene expression time series. $P$-values for these
test statistics are subsequently calculated, and then multiple testing
is performed based on these $p$-values while controlling the false
discovery rate (FDR) at a prespecified level, so that periodically
expressed genes are identified. Further improvements of this method
appear in \cite{Chen05,ALPHY05,GCM06,ALY07} and references
therein.

Global testing for periodicity dates back to the well-known Fisher's
test \cite{Fisher29}, which is based on the maximum value of the
normalized standard periodogram of the observed series. This test
enjoys some optimality properties as long as there is only one sinusoid
under the alternative and its frequency lies on the Fourier grid $(0,
\frac{2\pi}{N}, \frac{4\pi}{N}, \ldots, \frac{2(N-1)\pi}{N})$.
Since then substantial extensions and improvements have been made in
the literature. Based on an adaptive set of largest normalized standard
periodogram values, an extension of Fisher's test was proposed in \cite
{Siegel80}. It is empirically shown to be more powerful than Fisher's
test when there are multiple periodicities with Fourier frequencies
under the alternative. In \cite{Chiu89}, a modified Fisher's test, in
which the maximum periodogram is normalized by a trimmed mean of the
periodograms (as suggested in \cite{Bolviken83}), was proposed and
analyzed. This test is also more powerful than Fisher's test when there
are multiple periodicities under the alternative. In \cite
{Davies1987}, a test statistic was proposed against the alternative
when there is a single sinusoid whose frequency is unknown and not
necessarily on the Fourier grid. A more general hypothesis test is
given in \cite{JN15}, where the signal of interest consists of a fixed
number of sinusoids. The higher criticism test proposed in \cite
{DJ2004} can also be applied to the standard periodogram for global
detection, and it enjoys certain asymptotic minimaxity properties
against alternatives in which the periodicities are sparse and consist
of Fourier frequencies.

\setcounter{footnote}{1}

In the existing literature on global testing for periodicity, either
the contributing sinusoids are constrained to have Fourier frequencies,
or the number of periodicities is fixed and small compared to the
sample size. The goal of the present paper is to construct a new test
based on an \emph{over-sampled periodogram} that adapts to \emph{a
growing number of} sinusoids with \emph{general} frequencies. The
focus of our work is to establish the asymptotic optimality of our method.

%

\subsection{Methodology}
\label{semethod}

Our discussion throughout the paper is focused on  complex-valued time
series. As indicated in Section~4.1 of \cite{BD91} and Section~1.5 of
\cite{Fuller95}, complex-valued time series are sometimes more
convenient for analysis. Moreover, complex-valued or bivariate time
series arise naturally in modern data analysis such as functional MRI,
blood-flow and oceanography; see, for example, \cite{RW08,SOLE13} and
the references therein. Periodicity detection in real-valued series
will be briefly discussed in Section~\ref{secdiscussion}.


For ease of analysis, we slightly simplify the complex time series
model~\eqref{eqcomplexmodel} as follows:
%
\begin{equation}
\label{eqmeasurement} \mathbf{y}=\mathbf{X}\bolds{\beta}+\mathbf{z},
\end{equation}
where the design matrix $\mathbf{X} \in\mathbb{C}^{N \times p}$ with $p
\gg N$ is an \emph{extended discrete Fourier transform} (EDFT) matrix,
that is,
%
\begin{equation}
\label{eqextendedDFT} X_{jk}:=\frac{1}{\sqrt{p}}e^{-\sfrac{2 \pi i (k-1)(j-1)}{p}}, \qquad j=1,
\ldots, N,k=1, \ldots, p.
\end{equation}
The vector of coefficients $\bolds{\beta} \in\mathbb{C}^p$ contains
information of magnitudes and phases in \eqref{eqcomplexmodel}, and
its sparsity under the alternative, denoted as $s$, is assumed to be
\emph{unknown}. The noise level $\sigma$ is assumed to be known, and
we let $\sigma=1$ throughout the paper without loss of generality.

A distinct feature of our model in \eqref{eqmeasurement} is that the
value of $p$ can be arbitrarily large, and is assumed to be \emph
{unknown}. This implies that the design matrix $\mathbf{X}$ is actually
unavailable. Moreover, adjacent columns in $\mathbf{X}$ are nearly
parallel, which is different from the common assumption in
high-dimensional regression models in which the columns of the design
matrices are pairwise incoherent. A broad class of combinations of
periodicities can be represented by the mean $\mathbf{X}\bolds{\beta}$.
When $\bolds{\beta}$ is $s$-sparse, $\mathbf{X}\bolds{\beta}$ is a
superposition of $s$ complex sinusoids. The global test for periodicity
is therefore modeled as
%
\begin{equation}
\label{eqglobaltesting}
H_0: \bolds{\beta}=\mathbf{0} \quad\mbox{versus}\quad
H_1: \bolds {\beta } \neq\mathbf{0} \quad\mbox{and}\quad\mbox{$\bolds{\beta}$ is
sparse}.
\end{equation}

We now define the over-sampled periodogram for complex-valued time
series, which turns out to be surprisingly simple. For some integer
$q$, define
%
\begin{equation}
\label{eqmatrixU} \mathbf{U}= \biggl(\frac{1}{\sqrt{N}}e^{2\pi i \sfrac
{(m-1)(j-1)}{q}}
\biggr)_{1\leq m\leq q, 1\leq j\leq N}
\end{equation}
whose row vectors are normalized in $2$-norm, and set
%
\begin{equation}
\label{eqoversampledperiodogram} \mathbf{v}=\mathbf{U}\mathbf{y} \quad\mbox{and}\quad
I_m=|v_m|^2, \qquad m=1, \ldots, q.
\end{equation}
By letting $q>N$, $\{I_1, \ldots, I_q\}$ is an over-sampled
periodogram. Define
%
\begin{equation}
\label{eqHCt} \mathrm{HC}(t):= \frac{ \sum_{m=1}^{q} 1_{\{\sqrt{I_m} \geq t\}} -
q\bar{\Psi}(t)}{\sqrt{ q\bar{\Psi}(t)(1- \bar{\Psi}(t))}},
\end{equation}
where $\bar{\Psi}(t)=\mathbb{P}(|z|\geq t)=e^{-t^2}$ is the tail probability
of the standard complex normal variable as shown in Lemma~\ref
{teocomplextail}, and $I_m$ is defined in \eqref
{eqoversampledperiodogram}. Notice that under the null, all $v_m$ are
standard complex-valued normal variables, which implies $\mathbb{E}\mathrm{HC}(t)=0$
for any fixed $t$. The proposed test statistic is defined as
%
\begin{equation}
\label{eqHCstar} \mathrm{HC}^*= \sup_{a \leq t \leq b} \mathrm{HC}(t),
\end{equation}
for which appropriate choices of the interval $[a, b]$ are discussed in
Sections~\ref{sectheory} and~\ref{secnumerics}. We will fix a
threshold level $T$, and reject $H_0$ if and only if $\mathrm{HC}^*>T$. This
test is referred to as the \emph{over-sampled periodogram higher
criticism} (OPHC) test.

An important question is how to choose the over-sampling rate $q/N$.
Let $\bolds{\theta}=\mathbb{E}\mathbf{v} = \mathbf{U} \mathbf
{X}\bolds{\beta}$.
Roughly speaking, the success of detection by the higher criticism
based on the sequence $\mathbf{v}$ depends on whether $\bolds{\theta}$
has $s$ nonzero elements with sufficiently large magnitudes. If the
frequencies are on the Fourier grid, the spikiness of $\bolds{\theta}$
is implied by the spikiness of $\bolds{\beta}$. For example, if $s=1$ and
\[
y_t = \frac{A}{\sqrt{N}} e^{- i (\sfrac{2\pi}{N}) t} + z_t,
\]
by choosing $q=N$, we can calculate that $\bolds{\theta}$ has sparsity
one, and $\|\bolds{\theta}\|_\infty= A$.
Therefore, the proposed test is desirable as long as $A$ is
sufficiently large. However, if the frequencies are off the Fourier
grid, then for $q=N$, the spikiness of $\bolds{\beta}$ may not imply
the spikiness of $\bolds{\theta}$. For example, if
\[
y_t = \frac{A}{\sqrt{N}}e^{- i (\sfrac{\pi}{N}) t} + z_t,
\]
and one chooses $q=N$, simple calculation yields
\[
\limsup_{N \rightarrow\infty} \max_{1 \leq m \leq N} |
\theta_m| \leq\frac{2}{\pi} A.
\]
This means the resulting $\bolds{\theta}$ is not as spiky as in the
case where the frequencies are on the Fourier grid, and then the
performance of higher criticism based on $\mathbf{v}$ may be not optimal.

In order to increase the spikiness of $\bolds{\theta}$, we propose to
choose the over-sampling rate $q/N=O(\log N)$. Our main result Theorem~\ref{upperboundsparse} guarantees that as long as the frequencies of
the complex sinusoids in the mean of $\mathbf{y}$ obey some minimum
separation condition, this over-sampling rate leads to an
asymptotically optimal global test. A key step in the proof is to show
that $\bolds{\theta}$ has $s$ significant nonzero
components.\footnote
{This is indicated in equation \eqref{eqpeaks}.} In other words, the
spikiness of $\bolds{\beta}$ is translated to the spikiness of
$\bolds{\theta}$. We emphasize that this over-sampling rate is
independent of
the grid parameter $p$ and the sparsity $s$.


The higher criticism method was originally coined by John Tukey and
introduced in Donoho and Jin \cite{DJ2004} for signal detection under
a sparse homoscedastic Gaussian mixture model, which was previously
studied in Ingster \cite{Ingster1999}. Cai, Jin and Low \cite{CJL07}
investigated minimax estimation of the nonnull proportion $\epsilon_n$
under the same model. Hall and Jin \cite{HJ2010} proposed a modified
version of the high criticism for detection with correlated noise with
known covariance matrices.
Cai, Jeng and Jin \cite{CJJ11} considered heteroscedastic Gaussian
mixture model and showed that the optimal detection boundary can be
achieved by a double-sided version of the higher criticism test. The
papers \cite{CDH05,CCHZ08} considered a related problem of detecting a
signal with a known geometric shape in Gaussian noise. Cai and Wu \cite
{caiwu2014} studied the detection of sparse mixtures in the setting
where the null distribution is known, but not necessarily Gaussian and
established the adaptive optimality of the higher criticism for the
detection of such general sparse mixtures.

In the special case in which $p=N$, that is, the frequencies are on the
grid, the design matrix becomes the orthogonal DFT matrix. Multiplying
the measurement by the inverse DFT matrix, the design matrix is reduced
to the identity design. Therefore, the problem becomes equivalent to
the standard sparse detection model discussed in \cite{Ingster1999,DJ2004}, and the standard higher criticism test proposed in \cite
{DJ2004} can be directly applied. Notice that in the OPHC test defined
above, choosing $q=N$ in \eqref{eqmatrixU} is equivalent to
multiplying the measurement by the inverse DFT, so there is no need to
over-sample the periodogram.

\subsection{Relation with global testing in linear models}
\label{secliterature}
If the dimension $p$ in \eqref{eqmeasurement} were known, the
hypothesis testing model \eqref{eqmeasurement} considered in the
present paper is also closely related to the global testing problem
under a linear model with sparse alternatives. It is helpful to review
some well-known results for the real-valued case in this line of research.

Consider the linear model: $\mathbf{y}=\mathbf{X}\bolds{\beta}+\bolds
{\epsilon}$, where $\mathbf{X} \in\mathbb{R}^{N \times p}$, $\bolds
{\beta} \in\mathbb{R}^p$ are the design matrix and regression
coefficients, respectively. The noise vector $\bolds{\epsilon} \in
\mathbb{R}^N$ is assumed to be i.i.d. Gaussian variables with mean $0$
and variance $1$. The global detection of $\bolds{\beta}$ is still
captured by the hypothesis test \eqref{eqglobaltesting}. In the
recently developed literature of high-dimensional statistics, $p$ is
comparable or much greater than $N$, while the parameter vector $\bolds
{\beta}$ is assumed to be sparse: $\|\bolds{\beta}\|_0=s \ll N$. The
tradeoff between the strength of the nonzero regression coefficients
and the sparsity, by which the detectability of $\bolds{\beta}$ is
captured, has been intensively studied in the literature.

In order to simplify the analysis, it is convenient to assume that the
nonzero components of $\bolds{\beta}$ have the same magnitude $A$. The
tradeoff between the signal strength and sparsity is reduced to a
quantitative relationship between $A$ and $s$ for fixed $N$ and $p$.
This relationship also depends closely on the properties of the design
matrix $\mathbf{X}$. There are two well-studied examples in the literature:
\begin{itemize}
\item Identity design matrix. When $N=p$ and $\mathbf{X}=\mathbf{I}$, the
detection boundary is given in \cite{Ingster1999,DJ2004}. Let
$A=\sqrt{2r\log p}$ and $s=p^{1-\alpha}$ with $\alpha\in[\frac
{1}{2}, 1]$, a higher criticism test is asymptotically powerful as long
as $r > \rho^*(\alpha)$, where the detection boundary function $\rho
^*$ is defined as:
%
\begin{equation}
\label{eqrhostar} \rho^*(\alpha)= %
\cases{\displaystyle (1 - \sqrt{1 -
\alpha})^2, & $\displaystyle\quad\alpha\in\bigl[\tfrac{3}{4}, 1 \bigr)$,
\vspace*{3pt}\cr
\displaystyle \alpha-\tfrac{1}{2},  & $\quad\alpha\in \bigl(\tfrac{1}{2},
\tfrac{3}{4} \bigr)$.}
\end{equation}
On the other hand, if $r < \rho^*(\alpha)$, all sequences of testing
procedures are asymptotically powerless, and thus the signal is not
detectable. The condition $\alpha>\frac{1}{2}$ is crucial. Otherwise,
the detectability of nonzero $\bolds{\beta}$ is not characterized by
the scaling $A=\sqrt{2r\log p}$.

\item Gaussian design matrix. Another carefully studied class of design
matrices are the Gaussian designs; that is, $\mathbf{X} \in\mathbb
{R}^{N \times p}$ has i.i.d. zero-mean normal variables with variance
$\frac{1}{p}$. This model appears in \cite{ITV10} and \cite
{ACP2011}. By denoting $A=\sqrt{\frac{2rp \log p}{N}}$ and $s=p^{1 -
\alpha}$, the detection boundary established in \cite{ITV10} is still
$r=\rho^*(\alpha)$ as in the case of identity design, provided $p^{1
- \alpha}\log(p) = o(\sqrt{N})$. A similar result is established in
\cite{ACP2011}.
\end{itemize}

For ease of presentation, we assume $N=p^{1 - \gamma}$ with $0 \leq
\gamma< 1$ throughout the paper. Then in the case of Gaussian designs,
the detection boundary $r=\rho^*(\alpha)$ holds when $(1+ \gamma)/2
< \alpha<1$. However, there is an ``unnatural'' property of the
detection boundary $\rho^*$ in this case: When $\gamma>0$, $r
\rightarrow\rho^*(\frac{1+ \gamma}{2})>0$ as $\alpha\rightarrow
\frac{1+ \gamma}{2}$. In Section~\ref{sectheory}, with a similar
setup of $\gamma$, $\alpha$ and $r$, under the condition $\frac{1+
\gamma}{2} < \alpha<1$, a new detection boundary is developed for the
model \eqref{eqmeasurement} with EDFT designs, as long as the support
of $\bolds{\beta}$ satisfies a mild minimum separation condition. To be
specific, the new detection boundary is defined as
%
\begin{equation}
\label{eqrhogammastar}
\rho_\gamma^*(\alpha)= %
\cases{\displaystyle \bigl(\sqrt{1
- \gamma} - \sqrt{1 - \alpha}\bigr)^2, & $\displaystyle\quad\alpha\in \biggl[\frac{3}{4}+
\frac{\gamma}{4}, 1 \biggr)$,
\vspace*{3pt}\cr
\displaystyle\alpha-\frac{1}{2} - \frac{\gamma}{2}, &$\displaystyle\quad\alpha\in \biggl[
\frac{1 + \gamma}{2}, \frac{3}{4}+\frac{\gamma}{4} \biggr)$.
}
\end{equation}
It enjoys the property $r \rightarrow\rho_\gamma^*(\frac{1+ \gamma
}{2})=0$ as $\alpha\rightarrow\frac{1+ \gamma}{2}$.\hspace*{1pt} A detailed
comparison between the detection boundary of EDFT designs and that of
Gaussian designs is also provided in Section~\ref{sectheory}.

\subsection{Structure of the paper}\label{sec13}
The rest of the paper is organized as follows: In Section~\ref{sectheory}, we give theoretical results for the proposed method. An
explicit detection boundary $r= \rho^*_\gamma(\alpha)$ is
established under a mild minimum separation assumption on the
underlying frequencies, and the asymptotic optimality of OPHC is
established. In Section~\ref{secnumerics}, numerical simulations
illustrate the efficacy of our approach. In the implementation of OPHC,
we compare the performances between $q/N=O(\log N)$, $q=N$ and $q=p$. A
summary of our main contributions is given in Section~\ref{secdiscussion}, along with some future research directions. All the
proofs are deferred to Section~\ref{secproofs}.

\section{Theoretical results}
\label{sectheory}
In this section, we aim to establish a sharp tradeoff between the
magnitudes and number of the nonzero components in $\bolds{\beta}$,
such that the OPHC test can successfully reject the null hypothesis
when the alternative is true. Under the alternative, we assume
$\operatorname{supp}(\bolds{\beta})=\{\tau_1, \ldots, \tau_s\}$,
where $1\leq\tau_1 <
\cdots< \tau_s \leq p$. This implies that the nonzero components of
$\bolds{\beta}$ are $\beta_{\tau_1}, \ldots, \beta_{\tau_s}$. If
we denote $\bolds{\tau}=
[
\tau_1, \ldots, \tau_s
]
^T$ and
%
\begin{equation}
\label{eqbetatilde} \bolds{\tilde{\beta}}= %
[ \tilde{\beta}_1,
\ldots, \tilde{\beta}_s ] %
^T: = %
[
\beta_{\tau_1}, \ldots, \beta_{\tau_s} ] %
^T,
\end{equation}
then under the alternative, the $s$-sparse signal $\bolds{\beta}$ is
uniquely parameterized by $(\bolds{\tau}, \bolds{\tilde{\beta}})$. The
distribution of the measurement $\mathbf{y}$ under the alternative is
therefore parameterized by $\bolds{\tau}$ and $\bolds{\tilde{\beta}}$,
denoted as $\mathbb{P}_{(\bolds{\tau}, \bolds{\tilde{\beta}})}$.
Under\vspace*{2pt} the
null, $\mathbf{y}$ consists of standard complex normal variables, denoted
as $\mathbb{P}_0$.

As discussed in Section~\ref{secliterature}, throughout the paper,
let $N=p^{1-\gamma}$ with fixed $\gamma\in[0, 1)$, and $s=p^{1 -
\alpha}$ with $\frac{1+\gamma}{2} < \alpha<1$. This implies that
$s<N^{\sfrac{1}{2}}$, which is consistent with the assumption in \cite
{ITV10,ACP2011}.
%

\subsection{Minimum separation condition}
We assume the distances between the indices of the nonzero components
of $\bolds{\beta}$, that is, $\tau_1, \ldots, \tau_s$, satisfy the
following minimum separation condition:
%
\begin{equation}
\label{eqminimumseparation}\hspace*{8pt} \Delta(\bolds{\tau}):=\frac{1}{p}\min\bigl\{|
\tau_{l+1} - \tau_l|:l=1, \ldots, s, \tau_{s+1}:=
\tau_1+p\bigr\} \geq\frac{\log^2 N}{N}.
\end{equation}
A similar minimum separation condition appears in the literature of
super-resolution; see \cite{Donoho1992,CF2014}.

This spacing condition holds asymptotically if the support is assumed
to be random. Assume that $\tau_1\leq\cdots\leq\tau_s$ are the
order statistics of independent uniformly distributed random variables
$a_1, \ldots, a_s$ in $\{1, \ldots, p\}$. For any $a, b \in\{1,
\ldots, p\}$, define the distance
%
\begin{equation}
\label{eqdistance} d(a, b)= \min \bigl({|a-b|}/{p}, 1-{|a-b|}/{p}\bigr).
\end{equation}
It is evident that $\Delta(\bolds{\tau}) = \min_{1\leq l_1< l_2\leq
s}d(a_{l_1}, a_{l_2})$. For any fixed $l_1< l_2$, and any $p_0 \in[0,
1]$, it is easy to see $\mathbb{P}(d(a_{l_1}, a_{l_2})<p_0) \leq
2p_0+\frac
{1}{p}$. Since there are $\frac{s(s-1)}{2}$ pairs, we have
\begin{eqnarray*}
&&\mathbb{P} \Bigl(\min_{1\leq l_1< l_2\leq s}d(a_{l_1},
a_{l_2})<p_0 \Bigr)\\
&&\qquad\leq\sum_{1\leq l_1< l_2\leq s}
\mathbb{P} \bigl(d(a_{l_1}, a_{l_2})<p_0 \bigr)
\leq s(s-1) \biggl(p_0 + \frac{1}{2p}\biggr).
\end{eqnarray*}
By letting $p_0=\frac{\log^2 N}{N}$, we obtain $\mathbb{P}
(\Delta
(\bolds{\tau}) <\frac{\log^2 N}{N} )\leq\frac{s^2\log^2
N}{N} + \frac{s^2}{p}$. Recall that we assume the sparsity satisfies
$\frac{1+\gamma}{2} < \alpha<1$, where $N=p^{1 - \gamma}$ and
$s=p^{1 - \alpha}$. Then $\frac{s^2\log^2 N}{N} + \frac{s^2}{p}
\rightarrow0$ as $p \rightarrow\infty$. Therefore, \eqref
{eqminimumseparation} holds with probability tending to $1$. A simple
corollary is that with probability approaching $1$, all the indices
$a_1, \ldots, a_s$ are distinct.
%

\subsection{Detection boundary}
Recall that under the alternative, the distribution of the observation
$\mathbf{y}$ is parameterized by $(\bolds{\tau}, \bolds{\tilde
{\beta }})$. We assume that\footnote{As discussed in Section~\ref{secliterature}, it is assumed that $A=\sqrt{\frac{2 r p \log
p}{N}}$ in the literature of global detection boundaries under linear
models. The difference of $\sqrt{2}$ stems from the difference between
real-valued and complex-valued sequences.}
%
\begin{eqnarray}
 (\bolds{\tau}, \bolds{\tilde{\beta}}) &\in & \Gamma(p, N, s, r)
\nonumber
\\[-8pt]
\label{eqparameterspace}
\\[-8pt]
\nonumber
&:=&
\biggl\{ |\tilde{\beta}_1|=\cdots=|\tilde{\beta}_s|=A=
\sqrt{\frac{r p\log
p}{N}}, \Delta(\bolds{\tau}) \geq\frac{\log^2 N}{N} \biggr\}.
\end{eqnarray}
For the parameter space $\Gamma(p, N, s, r)$, we aim to establish a
new minimax detection boundary $r=\rho^*_\gamma(\alpha)$ defined as
in \eqref{eqrhogammastar}, when the sparsity level satisfies $\frac
{1+\gamma}{2} < \alpha<1$. Recall that the OPHC test defined by
\eqref{eqmatrixU}--\eqref{eqHCstar} is determined by the interval
$[a, b]$, the specific choice of $q=O(N \log N)$, and the threshold
$T$. In our theoretical analysis, we choose $q=N\lfloor\log N +1
\rfloor$, $[a, b]=[1, \sqrt{\log\frac{N}{3}}]$, and $T=\log^2 N$.
The OPHC test is therefore defined as
%
\begin{equation}
\Psi= I\bigl(\mathrm{HC}^*>\log^2 N\bigr).
\end{equation}
That is, the null hypothesis is rejected if and only if $\mathrm{HC}^*>\log^2
N$. This threshold is often too conservative in practice, and a more
reliable and useful threshold for finite samples can be chosen by Monte
Carlo simulations, which we will discuss in Section~\ref{secnumerics}.


%

Our first theorem is regarding the detectable region of $(\alpha, r)$,
in which the null can be successfully rejected asymptotically.
%
\begin{theorem}
\label{upperboundsparse}
In the measurement model \eqref{eqmeasurement}, suppose $N=p^{1 -
\gamma}$ with $\gamma\in[0, 1)$. Under the alternative, we assume
$s=p^{1 - \alpha}$ with $\frac{1+\gamma}{2}<\alpha<1$, and $(\tau,
\bolds{\tilde{\beta}})$ satisfies \eqref{eqparameterspace} with
parameter $r$. Suppose $\rho_\gamma^*$ is defined as in \eqref
{eqrhogammastar}. If $r>\rho^*_\gamma(\alpha)$, the OPHC test
defined by \eqref{eqmatrixU}--\eqref{eqHCstar} with $q=N\lfloor
\log N +1 \rfloor$ and $[a, b]= [1, \sqrt{\log\frac
{N}{3}} ]$ is asymptotically powerful:
\[
\lim_{N \rightarrow\infty} \Bigl(\mathbb{P}_0(H_0
\mbox{ is rejected})+\max_{(\bolds{\tau}, \bolds{\tilde{\beta}})\in\Gamma(p,
N, s, r)} \mathbb{P}_{(\bolds{\tau}, \bolds{\tilde{\beta}})}
(H_0 \mbox{ is accepted} ) \Bigr)=0.
\]
\end{theorem}

The most significant technical novelty in this paper lies in the proof
of Theorem~\ref{upperboundsparse}. In the analysis of $\mathrm{HC}^*$ under
the alternative, the mean and covariance structure of $\mathbf{v}$, which
is defined in \eqref{eqoversampledperiodogram}, requires more
careful calculation than in existing work, for example, \cite{ACP2011,HJ2010}. In particular, the estimation of $\mathbb{E}(v_1), \ldots,
\mathbb{E}(v_q)$
and the control of $\operatorname{Cov}(1_{|v_a|>t}, 1_{|v_b|>t})$ are treated
cautiously based on a variety of cases. In relevant calculations, the
structure of the EDFT design matrix $\mathbf{X}$ needs to be sufficiently
employed. Under the null, the $\mathrm{HC}^*$ statistic is related to the
standard $\mathrm{HC}^*$ statistic discussed in \cite{DJ2004}, so the analysis
is easier than the case of alternative.


The following theorem gives the lower bound for the testing problem.
%
\begin{theorem}
\label{lowerboundsparse}
Under the same setup of Theorem~\ref{upperboundsparse}, if $r<\rho
_\gamma^*(\alpha)$, then all sequences of hypothesis tests are
asymptotically powerless, that is,
\[
\lim_{N \rightarrow\infty} \Bigl(\mathbb{P}_0(H_0
\mbox{ is rejected})+\max_{(\bolds{\tau}, \bolds{\tilde{\beta}})\in\Gamma(p,
N, s, r)} \mathbb{P}_{(\bolds{\tau}, \bolds{\tilde{\beta}})}
(H_0 \mbox{ is accepted} ) \Bigr)=1.
\]
\end{theorem}

The proof of Theorem~\ref{lowerboundsparse} is relatively easy, and
it is given in the supplemental material. In fact, by taking advantage
of the specific structure of the EDFT matrix $\mathbf{X}$, the deduction
can be reduced to the case $\mathbf{X}=\mathbf{I}$. The classic lower bound
arguments in \cite{Ingster1999,DJ2004,HJ2010} can then be directly applied.

\begin{figure}

\includegraphics{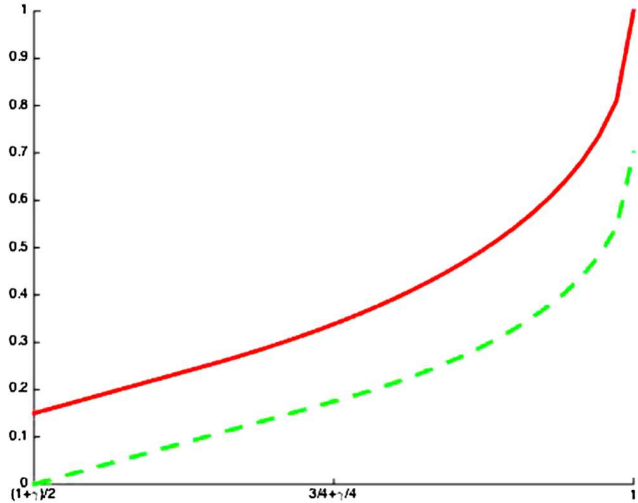}

\caption{Detection boundary functions $\rho^*(\alpha)$ (red solid
line) and $\rho^*_\gamma(\alpha)$ (green dashed line) for $\gamma
=0.3$ and $\frac{1+\gamma}{2}<\alpha<1$.}
\label{figdetectionboundaries}
\end{figure}

Theorems \ref{upperboundsparse} and \ref{lowerboundsparse}
together show that the proposed test is asymptotically optimal.
We now compare $\rho_\gamma^*$ with the detection boundary $\rho^*$
associated with the Gaussian designs established in \cite{ACP2011}. As
indicated in Section~\ref{sec13}, after normalizing the rows of the Gaussian
design, the magnitude parameter is denoted as $A=\sqrt{\frac{2rp \log
p}{N}}$. Notice that in our model the magnitude parameter is $A=\sqrt
{\frac{rp \log p}{N}}$, and the difference of $\sqrt{2}$ is due to
the distinction between real-valued and complex-valued models.
Therefore, it is fair to compare $\rho^*$ and $\rho_\gamma^*$
directly. It is obvious that $\rho_\gamma^*(\alpha)<\rho^*(\alpha
)$ for all $\frac{1+\gamma}{2}< \alpha<1$ as long as $\gamma>0$.
This implies that the detection boundary associated with the extended
DFT design matrix leads to milder trade-off between the rareness and
the weakness of the nonzero components of $\bolds{\beta}$ than that of
Gaussian designs. To illustrate their differences, the two detection
boundary functions are plotted in Figure~\ref{figdetectionboundaries} for $\gamma= 0.3$.


\section{Numerical simulations}
\label{secnumerics}
In this section, we study the empirical behavior of the OPHC test by
numerical simulations. In terms of computation, it is more convenient
to express the statistic as a function of $P_{(1)}\leq P_{(2)}\leq
\cdots\leq P_{(q)}$, which are the ordered $P$-values of $|v_1|, \ldots
, |v_q|$, that is, $P_{m} := \bar{\Psi}(|v_{m}|)=e^{-|v_{m}|^2}$. The
$\mathrm{HC}^*$ test in the following numerical simulations is defined as
%
\begin{equation}
\label{eqpractralOPHC} \mathrm{HC}^*=\max_{m: \sfrac{1}{q}\leq P_{(m)}<\sfrac{1}{2}} \frac{m
- qP_{(m)}}{\sqrt{qP_{(m)}(1 - P_{(m)})}},
\end{equation}
which is equivalent to choosing $[a, b]= [\sqrt{\log2}, \sqrt{\log
q}]$ in \eqref{eqHCstar}, instead of the theoretical choice $[a,
b]=[1, \sqrt{\log\frac{N}{3}}]$ defined in Section~\ref{sectheory}.


In the following, we compare the empirical testing powers of the OPHC
test with various choices of $q$.

\begin{figure}

\includegraphics{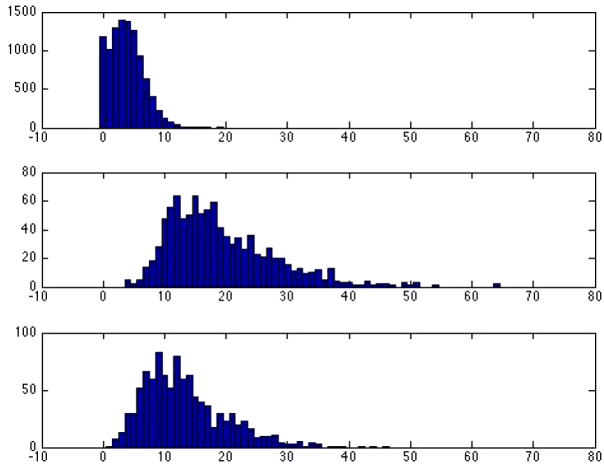}

\caption{Empirical distribution of the OPHC test statistic under the
null is plotted in the upper panel, where $N=1000$ and $q=14{,}000$. Under
the mixed alternative, we let $p=1{,}000{,}000$, $s=20$ and $r=0.3$.
Empirical distribution of the OPHC test statistic is plotted in the
middle panel with known $\sigma=1$, and in the lower panel with
estimated variance of noise.}
\label{figOPHC}
\end{figure}

First, let $N=1000$ and $q=2N\lfloor\log N + 1 \rfloor=14{,}000$. Then
the empirical distribution of the OPHC test statistic $\mathrm{HC}^*$ under the
null can be derived by Monte Carlo simulation with $1000$ independent
trials, which is shown in the upper panel of Figure~\ref{figOPHC}.

Under the alternative, we assume that $p=1{,}000{,}000$ and $s=20$. The
support of $\bolds{\beta}$ is distributed uniformly at random, and the
phase of the nonzero entries of $\bolds{\beta}$ are uniformly
distributed on $[0, 2\pi)$. All nonzero components of $\bolds{\beta}$
have the same magnitude $A=\sqrt{\frac{r p \log p}{N}}$ with $r=0.3$.

We first assume that the variance $\sigma^2=1$ is known. The resulting
empirical distribution of $\mathrm{HC}^*$ under the mixed alternative is plotted
in the middle panel of Figure~\ref{figOPHC} by $1000$ independent
trials. In $956$ trials of them, the empirical $P$-values are smaller
than $0.05$, by which the periodicities are successfully detected.

Let us discuss the case where the variance $\sigma^2=1$ is unknown.
Since it is necessary to make sure that the estimation of the variance
is consistent under the null, we use the mean square of $|y_t|$ as the
estimate. This estimator of $\sigma$ is adopted in order to make fair
numerical comparisons between different choices of $q$. Robust and
efficient estimation of $\sigma$ is an interesting problem. It has
been considered, for example, in \cite{dicker2014variance}, in the
case of Gaussian design. Efficient estimation of $\sigma$ in the
current setting is beyond the scope of this paper, and we leave it for
future research. The resulting empirical distribution of $\mathrm{HC}^*$ is
plotted in the lower panel of Figure~\ref{figOPHC}. Among the $1000$
trials, there are $745$ with empirical $P$-values smaller than $0.05$.

\begin{figure}[t]

\includegraphics{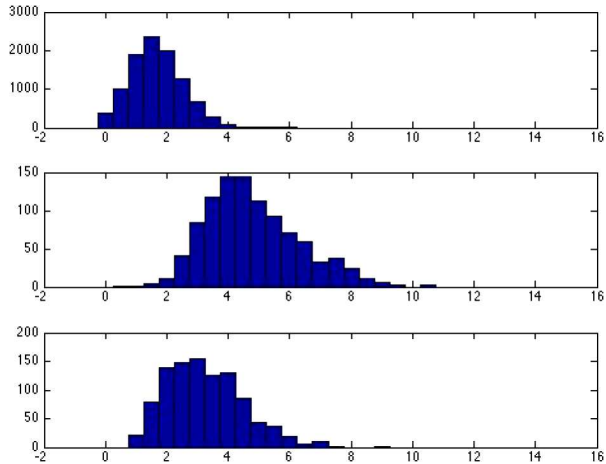}

\caption{Empirical distribution of the SPHC test statistic under the
null is plotted in the upper panel, where $N=1000$ and $q=1000$. Under
the mixed alternative, we let $p=1{,}000{,}000$, $s=20$ and $r=0.3$.
Empirical distribution of the SPHC test statistic is plotted in the
middle panel with known $\sigma=1$, and in the lower panel with
estimated variance of noise.}
\label{figPHCDJ}
\end{figure}

Next, we consider the OPHC test with $q=N=1000$. We refer to this test
as standard periodogram higher criticism (SPHC) test. The empirical
distribution of the SPHC test statistic under the null is plotted in
the upper panel of Figure~\ref{figPHCDJ}. The setup of the mixed
alternative is the same as in the experiments for the OPHC test
described before. Suppose the variance $\sigma^2=1$ is known. The
distribution of the SPHC test statistic under the alternative is
plotted in the middle panel of Figure~\ref{figPHCDJ} by $1000$
trials. In $867$ trials of them, the empirical $P$-values are smaller
than $0.05$. This is worse than the OPHC method with $q=14{,}000$, where
the periodicities are successfully detected in $956$ trials. When the
variance of the noise is unknown, we still estimate it by the mean
square of $|y_t|$, such that the estimate is consistent under the null.
The resulting empirical distribution of the SPHC test statistic is
plotted in the lower panel of Figure~\ref{figPHCDJ} based on $1000$
independent trials. In only $478$ trials among them, the empirical
$P$-values are smaller than $0.05$. This is also worse than OPHC where
the periodicities are successfully detected in $745$ independent trials.
\begin{figure}[t]

\includegraphics{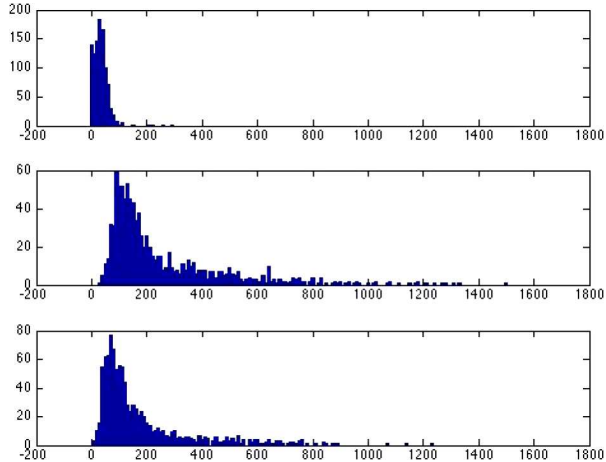}

\caption{Empirical distribution of the OPHC test statistic under the
null is plotted in the upper panel, where $N=1000$ and $q=1{,}000{,}000$.
Under the mixed alternative, we let $p=1{,}000{,}000$, $s=20$ and $r=0.3$.
Empirical distribution of this test statistic is plotted in the middle
panel with known $\sigma=1$, and in the lower panel with estimated
variance of noise.}
\label{figOPHCACP}
\end{figure}

Finally, we assume $p$ were known and consider the OPHC test with
$q=p=1{,}000{,}000$. In this case the OPHC test coincides with the method
proposed in \cite{ACP2011}. The empirical distribution of this test
statistic under the null by $1000$ independent trials is plotted in
the upper panel of Figure~\ref{figOPHCACP}. The setup of the mixed
alternative is the same as before. When the variance $\sigma^2=1$ is
known, the distribution of this test statistic under the alternative is
plotted in the middle panel of Figure~\ref{figOPHCACP} by $1000$
trials. In $949$ trials of them, the empirical $P$-values are smaller
than $0.05$, which is slightly worse than the OPHC method with
$q=14{,}000$ as mentioned before. When the variance of the noise is
unknown, with its estimation by the mean square of $|y_t|$, the
resulting empirical distribution of this test statistic is plotted in
the lower panel of Figure~\ref{figOPHCACP} based on $1000$
independent trials. In $741$ trials among them, the empirical $P$-values
are smaller than $0.05$, which is also slightly worse than the OPHC
method with $q=14{,}000$. Since $p$ is actually unknown, it would be more
convenient to choose $q=O(N\log N)$.

\section{Discussion}
\label{secdiscussion}

Motivated by periodicity detection in complex-valued time series
analysis, we investigated the hypothesis testing problem \eqref
{eqglobaltesting} under the linear model \eqref{eqmeasurement},
where the frequencies of the hidden periodicities are not necessarily
on the Fourier grid, and the number of sinusoids grows in $N$. The OPHC
test, a higher criticism test applied to the periodogram over-sampled
by $O(\log N)$, is proposed to solve this problem. In terms of theory,
by assuming that the frequencies satisfy a minimum separation
condition, a detection boundary between the rareness and weakness of
the sinusoids is explicitly established. Perhaps surprisingly, the
detectable region for the EDFT design matrix is broader than that of
Gaussian design matrices. For ease of exposition, we assume that $p$ is
finite but unknown, although $p$ is allowed to be infinity by slightly
modifying our argument.

Numerical simulations validate the choice $q=O(N \log N)$ by being
compared to the choice $q=N$, that is, the standard periodogram higher
criticism, and $q=p$, that is, excessively over-sampled periodogram
methods. In a recently published paper \cite{LS15}, it is shown that
higher criticism statistics might not be as powerful as Berk--Jones
statistics empirically. We find it interesting to investigate
alternative global testing methods for periodicity detection both
theoretically and empirically, but we leave this as future work.

The hypothesis testing problem considered in this paper is related to a
number of other interesting problems. We briefly discuss them here,
along with several directions for future research.

\subsection*{Statistical estimation of a large number of frequencies}
A related important statistical problem is to estimate the frequencies
of the periodicities in a given series. Sinusoidal regression methods
date back to 1795 by Prony \cite{Prony1795} with many later
developments including \cite{Schmidt1986,Cadzow1988,RK1989,HS1990},
to name a few. In the case where the number of frequencies is fixed and
few, numerous statistical analyses for frequency estimation have been
performed in the literature. For example, an insightful threshold
behavior of MLE was presented in \cite{QK94}. To determine the number
of frequencies, model selection methods are usually applied; see, for
example, \cite{Djuric1996,NK2011}. An extensive study on this subject
can be found in the classical text book \cite{QH01} and references
therein. In contrast, when the number of frequencies is large, although
sparse recovery \cite{CD1998,HSZF2012,FL2012,DB2013} and
total-variation minimization \cite{CF2014,CF2013} can be used for
frequency retrieval, their statistical efficiency is not clear. It is
interesting to develop both computationally and statistically efficient
methods to estimate a number of frequencies hidden in the observed sequence.


\subsection*{Sinusoidal denoising}
Compared to frequency estimation, denoising, that is, estimation of the
mean $\mathbf{X}\bolds{\beta}$ in \eqref{eqmeasurement}, is a
conceptually easier statistical task. In the recent papers \cite
{BTR2012,TBR2013}, SDP methods were shown to enjoy nearly-optimal
statistical properties. It would be interesting to establish
theoretically optimal methods.

\subsection{Testing for periodicity in real-valued series} Considering
the hypothesis test problem \eqref{eqtesting} in the real case \eqref
{eqrealmodel}, the OPHC test can be applied to the real sequence
$u_t$ by the idea of complexification, that is, transforming $u_t$ into
$y_t=u_t+iu_{t+n}$ for $t=1, \ldots, n$. Consequently, the mean of
$y_t$ amounts to a superposition of complex sinusoids, and the noise
part in $y_t$ consists of a sequence of complex white noise. Then the
hypothesis test problem is reduced to the complex case. It would be
interesting to investigate whether OPHC method can be applied to the
real series directly with potential statistical advantage.

\section{Proofs}
\label{secproofs}
This section is dedicated to the proofs of Theorems \ref
{upperboundsparse} and \ref{lowerboundsparse}. We begin by
collecting a few technical tools that will be used in the proof of the
main results.

\subsection{Preliminaries}
First, we formally introduce the concept of complex-valued multivariate
normal distribution.
%
\begin{definition}
\label{defcomplexnormal}
We say that $\mathbf{z}=\mathbf{x} + i\mathbf{y} \in\mathbb{C}^n$ is an
$n$-dimensional\break complex-valued multivariate normal vector with
distribution $\mathcal{CN}(\bolds{\mu}, \bolds{\Gamma}, \bolds
{\Omega })$, if $
\bigl[ {{\mathbf{x}} \atop {\mathbf{y}}}\bigr]
$ is an $2n$-dimensional real-valued normal vector, and $\mathbf{z}$ satisfies
\[
\mathbb{E}\mathbf{z} = \bolds{\mu}, \qquad\mathbb{E}\bigl((\mathbf{z}-\bolds {\mu})
(\mathbf{z}-\bolds{\mu})^*\bigr)=\bolds{\Gamma}, \qquad\mathbb {E}\bigl((\mathbf{z}-
\bolds{\mu}) (\mathbf{z}-\bolds{\mu})^T\bigr)=\bolds {\Omega}.
\]
Here, $\mathbf{X}^*$ denotes the complex conjugate transpose, while
$\mathbf{X}^T$ denotes the ordinary transpose. Moreover, we say a
complex-valued multivariate normal vector $\mathbf{z}$ is standard, if
\[
\mathbf{z} \sim\mathcal{CN}(\mathbf{0}, \mathbf{I}_n, \mathbf{0}).
\]
\end{definition}
The following lemma gives the tail distribution of the standard complex
normal variable, which turns out to be much neater than that of real
normal variables. Its proof is given in \cite{CEL2015}.

\begin{lemma}
\label{teocomplextail}
Suppose $z \sim\mathcal{C}\mathcal{N}(0, 1, 0)$ is the standard
circularly-symmetric complex normal variable. Then for any $t\geq0$,
$\bar{\Psi}(t):=\mathbb{P}(|z|\geq t)= e^{-t^2}$. In addition, if
$\mu\in
\mathbb{C}$ is fixed, we have
%
\begin{equation}
\label{eqnoncentrality} \frac{C_0}{1+(t-|\mu|)_+} e^{-(t - |\mu|)_+^2} \leq\mathbb{P}\bigl(|\mu + z|>t\bigr)
\leq e^{-(t - |\mu|)_+^2}
\end{equation}
for some positive numerical constant $C_0$. Here, $x_+:=\max(x, 0)$
for $x \in\mathbb{R}$, and $x_+^2$ is short for $(x_+)^2$.
\end{lemma}

Under the alternative hypothesis, we need an upper bound of the
variance of $\mathrm{HC}(t)$ for fixed $t$, for which the following lemma will
be applied for several times. The argument is standard in the
literature of normal comparison inequalities; see, for example, \cite
{LLR83,LS02}.

\begin{lemma}
\label{thmindcorr}
Suppose $
\bigl[ {w_1 \atop w_2}
\bigr]
\sim\mathcal{CN}(\mathbf{0}, \bolds{\Gamma}, \mathbf{0})$ is a
$2$-dimensional complex normal vector, where $\bolds{\Gamma}=
\bigl[ {1 \atop \bar{\xi}}\enskip {\xi \atop 1}\bigr]
$.
\begin{longlist}[2.]
\item[1.]  For any fixed $a_1, a_2 \in\mathbb{C}$ and $t>0$, there holds
\[
\operatorname{Cov}(1_{\{|w_1 - a_1|>t\}}, 1_{\{|w_2 - a_2|>t\}}) \leq\min \bigl(e^{-(t-|a_1|)_+^2},
e^{-(t-|a_2|)_+^2} \bigr).
\]
\item[2.]  If $|\xi|\leq\frac{1}{2}$, we obtain
\begin{eqnarray*}
&&\!\operatorname{Cov}(1_{\{|w_1 - a_1|>t\}}, 1_{\{|w_2 - a_2|>t\}})
\\
&&\!\!\qquad\leq C_0|\xi| \exp \biggl(-\frac{(t - |a_1|)_+^2+ (t -
|a_2|)_+^2}{1+|\xi|} \biggr)
\bigl(1+\bigl(t - |a_1|\bigr)_+\bigr) \bigl(1+\bigl(t - |a_2|\bigr)_+\bigr)
\end{eqnarray*}
for some positive numerical constant $C_0$.
\end{longlist}
\end{lemma}

The proof of this lemma is given in the supplemental material \cite{CEL2015}.


\begin{lemma}
\label{thmVanzuijlen}
Suppose $X_1, \ldots, X_n$ are $n$ i.i.d. random variables uniformly
distributed on $[0, 1]$. Then
\[
\mathbb{P} \biggl[\sup_{\sfrac{3}{n} \leq\rho\leq1 -\sfrac{3}{n}} \frac{|\sum_{j=1}^n 1_{\{X_j \leq\rho\}}- n\rho|}{\log n\sqrt
{n\rho(1-\rho)}} >
C_1 \biggr]\leq\frac{1}{\log^2 n}
\]
provided $n>N_0$. Here, $C_1$ and $N_0$ are positive numerical constants.
\end{lemma}
\begin{pf}
This is a weak version of existing concentration inequalities in \cite
{VanZuijlen1978} for ratio type empirical processes; see also \cite
{Alexander1987,GK06}.
\end{pf}
Finally, there is a simple and useful result which we will use in the
proofs for several times.
%
\begin{lemma}
\label{lmmdistance}
For any $a, b \in[0, 1]$, define $d(a, b)=\min(|a-b|, 1- |a-b|)$.
We have $4d(a, b)\leq\llvert 1 - e^{2\pi i (a - b)}\rrvert  \leq2\pi
d(a, b)$.
\end{lemma}
\begin{pf}
These inequalities can be obtained by comparing the length of arcs and
chords in the unit circle.
\end{pf}

\subsection{Proof of Theorem \texorpdfstring{\protect\ref{upperboundsparse}}{2.1}}
Suppose $\epsilon>0$ is an underdetermined positive parameter, which
will be specified later in order to establish the detectable region.
Denote $L=\lfloor\log N +1 \rfloor$, and hence $q=NL$. By the
definition of $\mathbf{v}$ \eqref{eqoversampledperiodogram} and the
definition of $\bolds{\tilde{\beta}}$ \eqref{eqbetatilde},
we have
\[
y_j = \frac{1}{\sqrt{p}} \sum_{l=1}^s
e^{-\sfrac{2 \pi i (j-1) (\tau
_l - 1)}{p}}\tilde{\beta}_l + z_j,\qquad j=1,\ldots, N,
\]
and
%
\begin{eqnarray}
v_m &=& \frac{1}{\sqrt{N}}\sum_{j=1}^N
e^{\sfrac{2 \pi i
(m-1)(j-1)}{q}} y_j
\nonumber
\\
&=&\frac{1}{\sqrt{Np}}\sum_{j=1}^N
\sum_{l=1}^s e^{2 \pi i
(j-1) (\vfrac{m-1}{q} - \vfrac{\tau_l - 1}{p} )}\tilde {
\beta}_l
\nonumber
\\[-8pt]
\label{eqthetaw}
\\[-8pt]
\nonumber
&&{} + \frac{1}{\sqrt{N}}\sum_{j=1}^N
e^{\sfrac{2 \pi i
(m-1)(j-1)}{q}} z_j
\\
&:= &\theta_m +w_m, \qquad m=1,\ldots, q.
\nonumber
\end{eqnarray}
It is obvious that $\theta_m$ is deterministic and $(w_1, \ldots,
w_q)$ is a $q$-dimensional complex multivariate normal vector. First,
since $z_j \sim\mathcal{CN}(0, 1, 0)$ are independent, we know
$\mathbb{E}
z_j^2 = 0$ and $\mathbb{E}z_j \bar{z}_j =1$.
This implies $\mathbb{E}w_m^2 =0$ and $\mathbb{E}w_m \bar{w}_m =1$
and, therefore,
$w_m \sim\mathcal{CN}(0, 1, 0)$. For any $1\leq m_1, m_2 \leq q$ and
$m_1 \neq m_2$, straightforward calculation gives $\mathbb{E}w_{m_1}w_{m_2}
=0$ and
\[
\mathbb{E}(w_{m_1} \bar{w}_{m_2})=\frac{1}{N} \sum
_{j=1}^N e^{\sfrac
{2 \pi
i (m_1 - m_2)(j-1)}{q}} =
\frac{1 - e^{\sfrac{2 \pi i N (m_1 -
m_2)}{q}}}{N (1 - e^{\sfrac{2 \pi i (m_1 - m_2)}{q}} )}.
\]
This\vspace*{1pt} implies $
\bigl[ { w_{m_1} \atop w_{m_2}}
\bigr]
\sim\mathcal{CN}\lleft(\mathbf{0},
\bigl[ {1 \atop  \bar{\xi}} \enskip {\xi \atop 1}
\bigr]
, \mathbf{0} \rright)$, where $\xi= \frac{1 - e^{\sfrac{2 \pi i N (m_1 -
m_2)}{q}}}{N (1 - e^{\sfrac{2 \pi i (m_1 - m_2)}{q}} )}$.
For all $m_1 \neq m_2$, by Lemma~\ref{lmmdistance},
%
\begin{equation}
\label{eqcorrelation} |\xi| \leq\frac{2}{N\llvert 1 - e^{\sfrac{2 \pi i (m_1 -
m_2)}{q}}\rrvert }\leq\frac{1}{2Nd (\vfrac{m_1 - 1}{q}, \vfrac
{m_2 - 1}{q} )}.
\end{equation}
Furthermore, when $\frac{m_1 - m_2}{L}=\frac{N(m_1 - m_2)}{q}$ is an
integer, we have
%
\begin{equation}
\label{eqindependence} \xi= \frac{1 - e^{\sfrac{2 \pi i N (m_1 - m_2)}{q}}}{N (1 -
e^{\sfrac{2 \pi i (m_1 - m_2)}{q}} )}=0.
\end{equation}
This implies that $w_{m_1}$ and $w_{m_2}$ are independent.

\subsubsection{Lower bound of $\mathbb{E}\mathrm{HC}(t)$ under the alternative}
%
%
\begin{longlist}
\item[\textit{Step} 1.]
We choose $m_1, \ldots, m_s \in\{1, \ldots, q\}$ such that $\frac
{m_1-1}{q}, \frac{m_2-1}{q}, \ldots, \frac{m_s-1}{q}$ are closest to
$\frac{\tau_1-1}{p}, \ldots, \frac{\tau_s-1}{p}$ under the metric
$d$, respectively. This implies that for each $1\leq l\leq s$,
%
\begin{equation}
\label{eqml} d \biggl(\frac{m_l -1}{q} , \frac{\tau_l - 1}{p} \biggr) \leq
\frac
{1}{2q} \leq\frac{1}{2N\log N}.
\end{equation}
Since $\bolds{\tau}$ satisfies the minimum separation condition as
indicated in \eqref{eqparameterspace}, $m_1, \ldots,\break  m_s$ must be distinct.
\end{longlist}

Therefore, for $1 \leq\nu\leq s$,
\begin{eqnarray*}
\theta_{m_{\nu}} & =& \mathop{\sum_{1\leq l \leq s }}_{l \neq\nu}
\frac{\tilde{\beta}_l}{\sqrt{Np}} \sum_{j=1}^N
e^{2 \pi i
(j-1) (\vfrac{m_\nu-1}{q} - \vfrac{\tau_l - 1}{p} )}\\
&&{}+\frac
{\tilde{\beta}_{\nu}}{\sqrt{Np}} \sum_{j=1}^N
e^{2 \pi i
(j-1) (\vfrac{m_\nu-1}{q} - \vfrac{\tau_{\nu} - 1}{p} )}
\\
&:=& S_1+S_2.
\end{eqnarray*}
First, as to $S_1$, we have
\begin{eqnarray*}
|S_1|&=& \frac{1}{\sqrt{Np}}\Biggl\llvert \sum
_{1\leq l \leq s , l \neq\nu
} \tilde{\beta}_l \sum
_{j=1}^N e^{2 \pi i (j-1) (\vfrac{m_{\nu
}-1}{q} - \vfrac{\tau_l - 1}{p} )}\Biggr\rrvert
\\
&=& \frac{1}{\sqrt{Np}}\biggl\llvert \sum_{1\leq l \leq s, l \neq\nu}
\tilde{\beta}_l \frac{1- e^{ 2 \pi i N (\vfrac{m_{\nu}-1}{q} -
\vfrac{\tau_l - 1}{p} )}}{1- e^{ 2 \pi i  (\vfrac{m_{\nu
}-1}{q} - \vfrac{\tau_l - 1}{p} )}}\biggr\rrvert
\\
&\leq& \frac{\|\bolds{\beta}\|_\infty}{\sqrt{Np}} \sum_{1\leq l \leq s, l \neq\nu}
\frac{2}{\llvert 1- e^{i 2 \pi
(\vfrac{m_{\nu} - 1}{q} - \vfrac{\tau_l -1}{p} )}\rrvert }
\\
& \leq & \frac{\|\bolds{\beta}\|_\infty}{\sqrt{Np}} \sum_{1\leq l \leq s, l \neq\nu}
\frac{2}{4d (\vfrac{m_{\nu} -
1}{q}, \vfrac{\tau_l -1}{p} )}.
\end{eqnarray*}
The last inequality is due to Lemma~\ref{lmmdistance}. Let us now
bound\hspace*{-3pt}
\[
\sum_{1\leq l \leq s, l \neq\nu} \frac{2}{4d
(\vfrac{m_{\nu} - 1}{q}, \vfrac{\tau_l -1}{p} )}.
\]
Since\vspace*{-3pt}
\begin{eqnarray*}
&& d \biggl(\frac{m_{\nu} - 1}{q}, \frac{\tau_l -1}{p} \biggr)\\
&&\qquad\geq   d \biggl(
\frac{\tau_l-1}{p}, \frac{\tau_{\nu}-1}{p} \biggr) - d \biggl(\frac{m_{\nu} - 1}{q},
\frac{\tau_{\nu}-1}{p} \biggr)
\\
&&\qquad\geq  \min\bigl(|\nu- l|, s- |\nu- l|)\Delta(\bolds{\tau}\bigr) - \frac{1}{2q}
\\
&&\qquad \geq  \min\bigl(|\nu- l|, s - |\nu- l|\bigr) \biggl(\Delta(\bolds{\tau}) -
\frac{1}{2q} \biggr),
\end{eqnarray*}
we have
\begin{eqnarray*}
&&\sum_{1\leq l \leq s, l \neq\nu} \frac{2}{4d (\vfrac
{m_{\nu} - 1}{q}, \vfrac{\tau_l -1}{p} )} \\
&&\qquad\leq \Biggl(\sum
_{a=1}^{\lfloor\sfrac{s}{2} \rfloor}\frac{1}{a} \Biggr)
\frac
{1}{\Delta(\bolds{\tau}) - \afrac{1}{2q}} \leq\frac{\log s}{\Delta
(\bolds{\tau}) - \afrac{1}{2q}} .
\end{eqnarray*}
As a result $|S_1| \leq\frac{\|\bolds{\beta}\|_\infty}{\sqrt{Np}}
\frac{\log s}{\Delta(\bolds{\tau}) - \afrac{1}{2q}}$.
As to $S_2$, we have
\begin{eqnarray*}
\biggl\llvert S_2 - \frac{\sqrt{N}\tilde{\beta}_{\nu}}{\sqrt{p}}\biggr\rrvert &=& 
\frac{\tilde{\beta}_{\nu}}{\sqrt{Np}} \Biggl(\sum_{j=1}^N
e^{2 \pi i (j-1) (\vfrac{m_{\nu}-1}{q} - \vfrac{\tau
_{\nu} - 1}{p} )} - N \Biggr)\Biggr\rrvert
\\
& \leq & \frac{\|\bolds{\beta}\|_\infty}{\sqrt{Np}} \sum_{j=1}^N
\bigl\llvert e^{2 \pi i (j-1) (\vfrac{m_{\nu}-1}{q} - \vfrac{\tau
_{\nu} - 1}{p} )} - 1\bigr\rrvert
\\
& =& \frac{\|\bolds{\beta}\|_\infty}{\sqrt{Np}} \sum_{j=1}^N \bigl
\llvert e^{2 \pi i (j-1)d (\vfrac{m_{\nu}-1}{q} , \vfrac{\tau_{\nu} -
1}{p} )} - 1\bigr\rrvert .
\end{eqnarray*}
The last inequality is by the definition of $d$ in Lemma~\ref
{lmmdistance}. By \eqref{eqml}, we have $d (\frac{m_{\nu
}-1}{q} , \frac{\tau_{\nu} - 1}{p} ) \leq\frac{1}{2q}$. By
Lemma~\ref{lmmdistance}, there holds
\begin{eqnarray*}
\bigl\llvert e^{2 \pi i (j-1)d (\vfrac{m_{\nu}-1}{q} , \vfrac{\tau
_{\nu} - 1}{p} )} - 1\bigr\rrvert &\leq &  d \biggl((j-1)2\pi d \biggl(
\frac{m_{\nu}-1}{q} , \frac{\tau_{\nu} - 1}{p} \biggr), 0 \biggr) \\
&\leq & \frac{(j - 1)\pi}{q}.
\end{eqnarray*}
This implies that
\begin{eqnarray*}
&& \biggl\llvert S_2 - \frac{\sqrt{N}\tilde{\beta}_{\nu}}{\sqrt{p}}\biggr\rrvert \leq
\frac{\|\bolds{\beta}\|_\infty}{\sqrt{Np}}\frac{\pi
}{q}\frac
{N(N-1)}{2}.
\end{eqnarray*}
In summary,
\begin{eqnarray*}
\biggl\llvert \theta_{m_\nu} -\sqrt{\frac{N}{p}} \tilde{
\beta}_\nu \biggr\rrvert &\leq & \frac{\|\bolds{\beta}\|_\infty}{\sqrt{Np}} \biggl(
\frac{\log s}{\Delta(\bolds{\tau})- \afrac{1}{2q}} + \frac{\pi
N(N-1)}{2q} \biggr)
\\
&\leq & \frac{\sqrt{r \log p}}{N} \biggl(\frac{\log s}{\vfrac{\log^2
N}{N} - \afrac{1}{2 N\log N}} + \frac{\pi N(N-1)}{2N\log N} \biggr).
\end{eqnarray*}
Noticing $N=p^{1 - \gamma}$ and $s=p^{1 - \alpha}$, for any fixed
$\epsilon>0$, we have
%
\begin{equation}
\label{eqpeaks} \biggl\llvert \theta_{m_\nu} -\sqrt{\frac{N}{p}}
\tilde{\beta}_\nu \biggr\rrvert \leq\epsilon^2 \sqrt{r
\log p}, \qquad\nu=1, \ldots, s,
\end{equation}
provided $p > C(\gamma, \alpha, r, \epsilon)$, which is a constant
only depending on $\gamma$, $\alpha$, $r$ and~$\epsilon$.

\begin{longlist}
\item[\textit{Step} 2.]
Define $D_l \subset\{1, \ldots, q\}$ for $1 \leq l \leq s$ as
\[
D_l = \biggl\{ m : 1\leq m\leq q, d \biggl(\frac{m}{q},
\frac{\tau
_l}{p} \biggr)<\frac{\sqrt{\log N}}{N} \biggr\},
\]
and $D= \bigcup_{1\leq l \leq s} D_l$. Since $\Delta(\bolds
{\tau}) \geq\frac{\log^2 N}{N} > \frac{2\sqrt{\log N}}{N}$, $D_1,
\ldots, D_s$ are disjoint\vspace*{1pt} subsets of $\{1, \ldots, q\}$. For each
index $m \in D_\nu$, $\nu=1, \ldots, s$, since
\begin{eqnarray*}
\theta_{m} & =& \mathop{\sum_{1\leq l \leq s}}_{l \neq\nu}
\frac
{\tilde{\beta}_l}{\sqrt{Np}} \sum_{j=1}^N
e^{2 \pi i (j-1)
(\vfrac{m-1}{q} - \vfrac{\tau_l - 1}{p} )}\\
&&{}+\frac{\tilde{\beta
}_{\nu}}{\sqrt{Np}} \sum_{j=1}^N
e^{2 \pi i (j-1)
(\vfrac{m-1}{q} - \vfrac{\tau_{\nu} - 1}{p} )} ,
\end{eqnarray*}
we have
\begin{eqnarray*}
|\theta_m|& \leq & \Biggl\llvert \frac{\tilde{\beta}_\nu}{\sqrt{Np}} \sum
_{j=1}^N e^{2 \pi i (j -1)  (\vfrac{m-1}{q}- \vfrac{\tau
_\nu-1}{p} )}\Biggr\rrvert \\
&&{}+\Biggl
\llvert \mathop{\sum_{1\leq l \leq s}}_{l \neq\nu} \frac{\tilde{\beta}_l}{\sqrt{Np}} \sum
_{j=1}^N e^{2 \pi i (j - 1)  (\vfrac{m-1}{q} - \vfrac{\tau_l-1}{p}
)}\Biggr\rrvert
\\
&\leq & \sqrt{\frac{N}{p}}\|\bolds{\beta}\|_\infty+
\frac{\|\bolds
{\beta}\|_\infty}{\sqrt{Np}} \mathop{\sum_{1\leq l \leq s}}_{l
\neq\nu} \frac{\llvert 1- e^{2 \pi i N (\vfrac{m-1}{q} - \vfrac
{\tau_l-1}{p} )}\rrvert }{\llvert 1- e^{2 \pi i  (\vfrac
{m-1}{q} - \vfrac{\tau_l-1}{p} )}\rrvert }
\\
& \leq & \sqrt{\frac{N}{p}}\|\bolds{\beta}\|_\infty+
\frac{\|\bolds
{\beta}\|_\infty}{\sqrt{Np}} \mathop{\sum_{1\leq l \leq s }}_{l
\neq\nu} \frac{1}{2d (\vfrac{m-1}{q}, \vfrac{\tau_l -
1}{p} )}.
\end{eqnarray*}
%
%
Notice that
\begin{eqnarray*}
d \biggl(\frac{m-1}{q}, \frac{\tau_l-1}{p} \biggr) &\geq &  d \biggl(
\frac{\tau_\nu-1}{p}, \frac{\tau_l-1}{p} \biggr) - d \biggl(\frac
{m - 1}{q} -
\frac{\tau_\nu-1}{p} \biggr)
\\
&\geq & \min\bigl(|\nu- l|, r- |\nu- l|\bigr)\Delta(\bolds{\tau}) - \frac
{\sqrt{\log N}}{N}
\\
& \geq& \min\bigl(|\nu- l|, r- |\nu- l|\bigr) \biggl(\Delta(\bolds{\tau}) -
\frac{\sqrt{\log N}}{N} \biggr).
\end{eqnarray*}
This implies that\vspace*{-5pt}
\begin{eqnarray*}
&&\mathop{\sum_{1\leq l \leq s}}_{l \neq\nu} \frac{1}{2d
(\vfrac{m-1}{q}, \vfrac{\tau_l - 1}{p} )} \\[-3pt]
&&\qquad\leq \Biggl(\sum
_{a=1}^{\lfloor\sfrac{s}{2} \rfloor}\frac{1}{a} \Biggr)
\frac
{1}{\Delta(\bolds{\tau}) - \sfrac{\sqrt{\log N}}{N}} \leq\frac
{\log
s}{\Delta(\bolds{\tau}) - \sfrac{\sqrt{\log N}}{N}}.
\end{eqnarray*}
In summary,\vspace*{-3pt}
\begin{eqnarray*}
|\theta_m| &\leq & \sqrt{\frac{N}{p}}\|\bolds{\beta}
\|_\infty+ \frac
{\|\bolds{\beta}\|_\infty}{\sqrt{Np}} \frac{\log s}{\Delta(\bolds
{\tau}) - \sfrac{\sqrt{\log N}}{N}}
\\
&=& \sqrt{r \log p} \biggl(1+\frac{\log s}{N \Delta(\bolds{\tau}) -
\sqrt{\log N}} \biggr)
\\
& \leq & \sqrt{r \log p} \biggl(1+\frac{\log s}{\log^2 N - \sqrt{\log
N}} \biggr).
\end{eqnarray*}
Similarly, when $p> C(\gamma, \alpha, r, \epsilon)$, we have
%
\begin{equation}
\label{eqthetaD} \llvert \theta_{m}\rrvert \leq\sqrt{r \log p}
\bigl(1+\epsilon ^2 \bigr) \qquad\forall m \in D.
\end{equation}
\item[\textit{Step} 3.]
In this step, we aim to give a uniform upper bound of $\theta_m$ for
all $m \in D^c$. Straightforward calculation yields
\begin{eqnarray*}
\llvert \theta_m\rrvert &=& \Biggl\llvert \sum
_{1\leq l \leq s } \frac
{\tilde{\beta}_l}{\sqrt{Np}} \sum_{j=1}^N
e^{2 \pi i (j-1)
(\vfrac{m-1}{q} - \vfrac{\tau_l - 1}{p} )}\Biggr\rrvert
\\
&\leq & \frac{\|\bolds{\beta}\|_\infty}{\sqrt{Np}} \sum_{1\leq l \leq s}
\frac{\llvert 1- e^{i 2 \pi N (\vfrac{m-1}{q} -
\vfrac{\tau_l-1}{p} )}\rrvert }{\llvert 1- e^{i 2 \pi
(\vfrac{m-1}{q} - \vfrac{\tau_l-1}{p} )}\rrvert } \\
&\leq & \frac{\|
\bolds{\beta}\|_\infty}{\sqrt{Np}} \sum_{1\leq l \leq s}
\frac{1}{2d (\vfrac{m-1}{q}, \vfrac{\tau_l - 1}{p} )}.
\end{eqnarray*}
We now aim to bound $\sum_{1\leq l \leq s} \frac{1}{2d
(\vfrac{m-1}{q}, \vfrac{\tau_l - 1}{p} )}$. We consider the
position of $\frac{m - 1}{q}$ on $T=[0, 1]/\{0\sim1\}$ relative to
$\frac{\tau_1 -1}{p}, \ldots, \frac{\tau_s -1}{p}$. Suppose on
$T=[0, 1]/\{0\sim1\}$, $\frac{m-1}{q}$ is located between $\frac
{\tau_j -1}{p}$ and $\frac{\tau_{j+1} -1}{p}$ (recall that $\tau
_{s+1}=\tau_1$). Since $m \in D^c$, we have $d(\frac{m - 1}{q}, \frac
{\tau_j -1}{p})\geq\frac{\sqrt{\log N}}{N}$ and $d(\frac{m -
1}{q}, \frac{\tau_{j+1} - 1}{p})\geq\frac{\sqrt{\log N}}{N}$. The
next adjacent location parameters $\tau_{j-1}$ and $\tau_{j+2}$
satisfy $d(\frac{m - 1}{q},\break \frac{\tau_{j-1}-1}{p})\geq\Delta
(\bolds{\tau})+\frac{\sqrt{\log N}}{N}$ and $d(\frac{m - 1}{q},
\frac
{\tau_{j+2} - 1}{p})\geq\Delta(\bolds{\tau})+\frac{\sqrt{\log
N}}{N}$, etc. Then we have
\begin{eqnarray*}
\sum_{1\leq l \leq s} \frac{1}{2d (\vfrac{m-1}{q}, \vfrac
{\tau_l - 1}{p} )} &\leq & \sum
_{a=0}^{\lfloor\vfrac{s - 1}{2}
\rfloor}\frac{1}{a\Delta(\bolds{\tau}) + \sfrac{\sqrt{\log
N}}{N}}\\
&\leq & \frac{N}{\sqrt{\log N}}
+ \frac{\log s}{\Delta(\bolds
{\tau})}.
\end{eqnarray*}
In summary,
\[
|\theta_m| \leq\frac{\|\bolds{\beta}\|_\infty}{\sqrt{Np}} \biggl(\frac{N}{\sqrt{\log N}} +
\frac{\log s}{\Delta(\bolds{\tau
})} \biggr)=\sqrt{r \log p} \biggl(\frac{1}{\sqrt{\log N}} +
\frac
{\log s}{\log^2 N} \biggr).
\]
Similarly, when $p> C(\gamma, \alpha, r, \epsilon)$, we have
%
\begin{equation}
\label{eqthetanotinD} |\theta_m|\leq\epsilon^2 \sqrt{r \log
p} \qquad\forall m \in D^c.
\end{equation}

\item[\textit{Step} 4.]
We are now ready to derive a lower bound of $\mathbb{E}\mathrm{HC}(t)$. Recall that
%
\begin{equation}
\label{eqHCtt} \mathrm{HC}(t) = \frac{ \sum_{m=1}^q 1_{\{|v_m|>t\}} - q \bar{\Psi
}(t)}{\sqrt{ q \bar{\Psi}(t)(1- \bar{\Psi}(t))}}.
\end{equation}
By Lemma~\ref{teocomplextail}, we have
\begin{eqnarray*}
\mathbb{E}\mathrm{HC}(t) &\geq & \frac{\sum_{l=1}^s \mathbb{P}(|\theta_{m_l}
+ w_{m_l}|>t) -
s \bar{\Psi}(t)}{\sqrt{q \bar{\Psi}(t)(1- \bar{\Psi}(t))}}
\\
& \geq & \frac{\sum_{l=1}^s \afrac{C_0}{1+t}e^{-(t - |\theta
_{m_l}|)_+^2} - s e^{-t^2}}{\sqrt{q \bar{\Psi}(t)(1 - \bar{\Psi}(t))}}.
\end{eqnarray*}
By \eqref{eqpeaks}, we have
\[
\min_{1\leq l \leq s} |\theta_{m_l}| \geq\bigl( 1-\epsilon
^2\bigr)\sqrt{r\log p},
\]
which implies
\[
\mathbb{E}\mathrm{HC}(t)\geq\frac{ s (\afrac{C_0}{1+t}e^{-(t - (
1-\epsilon
^2)\sqrt{r\log p})_+^2} - e^{-t^2} )}{\sqrt{q \bar{\Psi}(t)(1
- \bar{\Psi}(t))}}.
\]
Letting $t=\sqrt{\mu\log p}$, by $s=p^{1 - \alpha}$, $N=p^{1 -
\gamma}$ and $q=N\lfloor\log N +1 \rfloor$, there holds
%
\begin{equation}
\label{eqexpectation} \mathbb{E}\mathrm{HC}(\sqrt{\mu\log p})\geq\frac{1}{\operatorname{polylog}(p)}p^{\sfrac
{1}{2} - \alpha+ \sfrac{\gamma}{2} + \sfrac{\mu}{2} - (\sqrt{\mu}
- (1- \epsilon^2)\sqrt{r})_+^2},
\end{equation}
where $\operatorname{polylog}(p)$ is a polynomial of $\log p$.
\end{longlist}

\subsubsection{Upper bound of $\operatorname{Var}(\mathrm{HC}(t))$ under the
alternative}
By \eqref{eqthetaD} and \eqref{eqthetanotinD}, we have
\[
\cases{ \displaystyle \max_{1 \leq m \leq q} |\theta_m|
\leq \bigl(1+\epsilon ^2 \bigr)\sqrt{r \log p},
\vspace*{3pt}\cr
\displaystyle  \max_{m \in D^c} |\theta_m|\leq\epsilon^2
\sqrt{ r \log p}.}
\]
By the definition of $\mathrm{HC}(t)$ as in \eqref{eqHCtt}, simple calculation yields
\[
\operatorname{Var}\mathrm{HC}(t) = \frac{1}{q \bar{\Psi}(t)(1 - \bar{\Psi
}(t))} \sum_{1
\leq a, b \leq q}
\operatorname{cov} (1_{\{|\theta_a + w_a|>t\}}, 1_{\{
|\theta_b + w_b|>t\}} ).
\]
By equation \eqref{eqcorrelation}, when $d ( \frac{a - 1}{q},
\frac{b - 1}{q} ) \geq\frac{1}{N}$, we have
\[
\bigl|\mathbb{E}(w_a\bar{w_b})\bigr|\leq\frac{1}{2Nd (\vfrac{a - 1}{q},
\vfrac{b
- 1}{q} )} \leq
\frac{1}{2}.
\]
Then Lemma~\ref{thmindcorr} implies
\[
\operatorname{cov} (1_{\{|\theta_a + w_a|>t\}}, 1_{\{|\theta_b +
w_b|>t\}} ) \leq\frac{C_0\exp (- \vfrac{(t - |\theta
_a|)_+^2 + (t - |\theta_b|)_+^2}{2} )(1+t)^2}{2Nd (\vfrac{a
- 1}{q}, \vfrac{b - 1}{q} )}.
\]
On the other hand, when $d ( \frac{a - 1}{q}, \frac{b -
1}{q} ) < \frac{1}{N}$, Lemma~\ref{thmindcorr} implies
\[
\operatorname{cov} (1_{\{|\theta_a + w_a|>t\}}, 1_{\{|\theta_b +
w_b|>t\}} ) \leq e^{-(t - |\theta_a|)_+^2}.
\]
Now we bound $\operatorname{Var}\mathrm{HC}(t)$ by controlling
\begin{eqnarray*}
S_1 (t) &=&  \frac{1}{q \bar{\Psi}(t)(1 - \bar{\Psi}(t))} \\
&&{}\times\sum_{a \in
D}
\sum_{d ( \vfrac{a - 1}{q}, \vfrac{b - 1}{q} ) < \sfrac
{1}{N}} \operatorname{cov} (1_{\{|\theta_a + w_a|>t\}},
1_{\{|\theta_b
+ w_b|>t\}} ),
\\
S_2 (t) &=&  \frac{1}{q \bar{\Psi}(t)(1 - \bar{\Psi}(t))}\\
&&{}\times \sum_{a \in
D}
\sum_{d ( \vfrac{a - 1}{q}, \vfrac{b - 1}{q} ) \geq\sfrac
{1}{N}} \operatorname{cov} (1_{\{|\theta_a + w_a|>t\}},
1_{\{|\theta_b
+ w_b|>t\}} ),
\\
S_3 (t) &=&  \frac{1}{q \bar{\Psi}(t)(1 - \bar{\Psi}(t))}\\
&&{}\times \sum_{a \in
D^c}
\mathop{\sum_{d ( \vfrac{a - 1}{q}, \vfrac{b - 1}{q} ) < \sfrac{1}{N}}}_{b \in D^c} \operatorname{cov} (1_{\{|\theta_a + w_a|>t\}},
1_{\{|\theta_b + w_b|>t\}} ),
\end{eqnarray*}
and
\begin{eqnarray*}
S_4 (t)&=& \frac{1}{q \bar{\Psi}(t)(1 - \bar{\Psi}(t))} \\
&&{}\times\sum_{a \in
D^c}
\mathop{\sum_{d ( \vfrac{a - 1}{q}, \vfrac{b - 1}{q} ) \geq
\sfrac{1}{N} }}_{b \in D^c} \operatorname{cov} (1_{\{|\theta_a +
w_a|>t\}},
1_{\{|\theta_b + w_b|>t\}} ).
\end{eqnarray*}
By the symmetry between $a$ and $b$, we have
\[
\operatorname{Var}\mathrm{HC}(t) \leq2\bigl(S_1(t) + S_2(t)
\bigr) + S_3(t) + S_4(t).
\]

\begin{longlist}
\item[\textit{Step} 1: \textit{Upper bound for} $S_1 (t)$.]
For fixed $a \in D$,
\begin{eqnarray*}
\sum_{d ( \vfrac{a - 1}{q}, \vfrac{b - 1}{q} ) < \sfrac
{1}{N}} \operatorname{cov} (1_{\{|\theta_a + w_a|>t\}},
1_{\{|\theta_b
+ w_b|>t\}} ) &\leq & \frac{2q}{N} e^{-(t - |\theta_a|)_+^2}
\\
&\leq & \frac{2q}{N}e^{- (t - \sqrt{r \log p}(1+\epsilon^2) )_+^2}.
\end{eqnarray*}
Notice that $|D| \leq\frac{2qs\sqrt{\log N}}{N}$, we get
\[
S_1(t)\leq\frac{\sfrac{2qs\sqrt{\log N}}{N}}{q \bar{\Psi}(t)(1 -
\bar{\Psi}(t))}\frac{2q}{N}e^{- (t - \sqrt{r \log
p}(1+\epsilon^2) )_+^2}
\]
which implies
\[
S_1(\sqrt{\mu\log p}) \leq\operatorname{polylog}(p) p^{-\alpha+ \gamma+
\mu-(\sqrt{\mu} - (1+\epsilon^2)\sqrt{r})_+^2}.
\]

\item[\textit{Step} 2: \textit{Upper bound for} $S_2 (t)$.]
For fixed $a \in D$,
\begin{eqnarray*}
&& \sum_{ d ( \vfrac{a - 1}{q}, \vfrac{b - 1}{q} )
\geq\sfrac{1}{N}} \operatorname{cov} (1_{\{|\theta_a + w_a|>t\}},
1_{\{
|\theta_b + w_b|>t\}} )
\\
&&\qquad \leq \sum_{ d ( \vfrac{a - 1}{q}, \vfrac{b - 1}{q} )
\geq\sfrac{1}{N}}\frac{C_0\exp (- \vfrac{(t - |\theta_a|)_+^2
+ (t - |\theta_b|)_+^2}{2} )(1+t)^2}{2Nd (\vfrac{a - 1}{q},
\vfrac{b - 1}{q} )}
\\
&&\qquad  \leq\frac{C_0(1+t)^2e^{- (t - \sqrt{r \log p}(1+\epsilon
^2) )_+^2}}{2N}\\
&&\qquad\quad{}\times \sum_{d ( \vfrac{a - 1}{q}, \vfrac{b -
1}{q} ) \geq\sfrac{1}{N}}
\frac{1}{d (\vfrac{a - 1}{q},
\vfrac{b - 1}{q} )}
\\
&&\qquad \leq\frac{C_0(1+t)^2e^{- (t - \sqrt{r \log p}(1+\epsilon
^2) )_+^2}}{2N} 2\sum_{l=1}^q
\frac{q}{l}
\\
&&\qquad \leq\frac{C_0q \log q(1+t)^2e^{- (t - \sqrt{r \log
p}(1+\epsilon^2) )_+^2} }{N}.
\end{eqnarray*}
Notice that $|D| \leq\frac{2qs\sqrt{\log N}}{N}$, we get
\[
S_2(t)\leq\frac{\sfrac{2qs\sqrt{\log N}}{N}}{q \bar{\Psi}(t)(1 -
\bar{\Psi}(t))}\frac{C_0q \log q(1+t)^2e^{- (t - \sqrt{r \log
p}(1+\epsilon^2) )_+^2} }{N},
\]
which implies
\[
S_2(\sqrt{\mu\log p}) \leq\operatorname{polylog}(p) p^{-\alpha+ \gamma+
\mu-(\sqrt{\mu} - (1+\epsilon^2)\sqrt{r})_+^2}.
\]

\item[\textit{Step} 3: \textit{Upper bound for} $S_3 (t)$.]
For fixed $a \in D^c$ and any $b \in\{1, \ldots, q\}$
\[
\operatorname{cov} (1_{\{|\theta_a + w_a|>t\}}, 1_{\{|\theta_b +
w_b|>t\}} ) \leq e^{-(t - |\theta_a|)_+^2} \leq
e^{-(t -
\epsilon^2 \sqrt{r\log p})_+^2}.
\]
Then
\begin{eqnarray*}
S_3(t)& \leq & \frac{1}{q \bar{\Psi}(t)(1 - \bar{\Psi}(t))}
\sum_{a
\in D^c}
\mathop{\sum_{d ( \vfrac{a - 1}{q}, \vfrac{b - 1}{q} ) <
\sfrac{1}{N}}}_{b \in D^c} e^{-(t - \epsilon^2 \sqrt{r\log p})_+^2}
\\
&=& \frac{q (\sfrac{2q}{N}) e^{-(t - \epsilon^2 \sqrt{r\log p})_+^2} }{q
\bar{\Psi}(t)(1 - \bar{\Psi}(t))},
\end{eqnarray*}
which implies\vspace*{-3pt}
\[
S_3(\sqrt{\mu\log p}) \leq\operatorname{polylog}(p) p^{\mu- (\sqrt{\mu}
- \epsilon^2 \sqrt{r})_+^2}.
\]

\item[\textit{Step} 4: \textit{Upper bound for} $S_4 (t)$.]
For fixed $a \in D^c$,
\begin{eqnarray*}
&&\mathop{\sum_{b \in D^c}}_{d ( \vfrac{a - 1}{q}, \vfrac{b -
1}{q} ) \geq\sfrac{1}{N}} \operatorname{cov} (1_{\{|\theta_a +
w_a|>t\}},
1_{\{|\theta_b + w_b|>t\}} )
\\
&&\qquad \leq\mathop{\sum_{b \in D^c}}_{d ( \vfrac{a - 1}{q}, \vfrac{b -
1}{q} ) \geq\sfrac{1}{N}}\frac{C_0\exp (- \vfrac{(t -
|\theta_a|)_+^2 + (t - |\theta_b|)_+^2}{2} )(1+t)^2}{2Nd
(\vfrac{a - 1}{q}, \vfrac{b - 1}{q} )}
\\
&&\qquad \leq\frac{C_0(1+t)^2e^{- (t - \epsilon^2\sqrt{r \log
p} )_+^2}}{2N} \\
&&\qquad\quad{}\times\mathop{\sum_{b \in D^c}}_{d ( \vfrac{a - 1}{q},
\vfrac{b - 1}{q} ) \geq\sfrac{1}{N}}
\frac{1}{d (\vfrac{a
- 1}{q}, \vfrac{b - 1}{q} )}
\\
&&\qquad \leq\frac{C_0(1+t)^2e^{- (t - \epsilon^2\sqrt{r \log
p} )_+^2}}{2N} 2\sum_{l=1}^q
\frac{q}{l}
\\
&&\qquad \leq\frac{C_0q \log q(1+t)^2e^{- (t - \epsilon^2\sqrt{r \log
p} )_+^2} }{N}.
\end{eqnarray*}
Therefore, $S_4(t)\leq\frac{q}{q \bar{\Psi}(t)(1 - \bar{\Psi
}(t))}\frac{C_0q \log q(1+t)^2e^{- (t - \epsilon^2\sqrt{r \log
p} )_+^2} }{N}$,
which implies
\[
S_4(\sqrt{\mu\log p}) \leq\operatorname{polylog}(p) p^{\mu- (\sqrt{\mu}
- \epsilon^2 \sqrt{r})_+^2}.
\]
\item[\textit{Step} 5: \textit{Summary}.]
In summary, there holds
%
\begin{eqnarray}
&&\operatorname{Var}\mathrm{HC}(\sqrt{\mu\log p})\nonumber\\
\label{eqvariance}
 &&\qquad \leq  2
\bigl(S_1(\sqrt{\mu\log p}) + S_2(\sqrt{\mu\log p})\bigr)
+ S_3(\sqrt{\mu\log p}) + S_4(\sqrt{\mu \log p})
\\
&&\qquad\leq  \operatorname{polylog}(p) \bigl(p^{-\alpha+ \gamma+ \mu- (\sqrt{\mu
} - (1+\epsilon^2)\sqrt{r})_+^2} + p^{\mu- (\sqrt{\mu} - \epsilon
^2\sqrt{r})_+^2} \bigr).\nonumber
\end{eqnarray}
\end{longlist}

\subsubsection{Detectable region under the alternative and the values
of \texorpdfstring{$\epsilon$}{$epsilon$} and \texorpdfstring{$\mu$}{$mu$}}
By Chebyshev's inequality, we have
\begin{eqnarray*}
&& \mathbb{P} \bigl(\mathrm{HC}(\sqrt{\mu\log p}) < \mathbb {E}\mathrm{HC}(\sqrt{\mu\log
p}) - (\log N) \bigl(\operatorname{Var} \bigl(\mathrm{HC}(\sqrt{\mu\log p}) \bigr)
\bigr)^{\sfrac{1}{2}} \bigr) \\
&&\qquad\leq\frac{1}{\log^2 N}.
\end{eqnarray*}
By \eqref{eqexpectation} and \eqref{eqvariance}, to guarantee
$\mathbb{P}
 (\mathrm{HC}(\sqrt{\mu\log p})\leq\log^2 N )\leq\frac{1}{\log
^2 N}$, it suffices to require that $p > C(\gamma, \alpha, r,
\epsilon, \mu)$ is sufficiently large, and
\[
\cases{\displaystyle \frac{1}{2} - \alpha+ \frac{\gamma}{2}
+ \frac{\mu}{2} - \bigl(\sqrt {\mu} - \bigl(1- \epsilon^2\bigr)
\sqrt{r}\bigr)_+^2\vspace*{3pt}\cr
\displaystyle\quad >\frac{1}{2} \bigl(-\alpha+ \gamma+ \mu-
\bigl(\sqrt{\mu} - \bigl(1+\epsilon^2\bigr)\sqrt{r}
\bigr)_+^2 \bigr), \vspace*{3pt}
\cr
\displaystyle\frac{1}{2} -
\alpha+ \frac{\gamma}{2} + \frac{\mu}{2} - \bigl(\sqrt {\mu} - \bigl(1-
\epsilon^2\bigr)\sqrt{r}\bigr)_+^2 > \frac{1}{2}
\bigl(\mu- \bigl(\sqrt{\mu} - \epsilon^2\sqrt{r}\bigr)_+^2
\bigr)>0, } %
\]
which amounts to
%
\begin{equation}
\label{eqconstraints} %
\cases{\displaystyle 1 - \alpha> 2\bigl(\sqrt{\mu} -
\bigl(1- \epsilon^2\bigr)\sqrt{r}\bigr)_+^2 - \bigl(\sqrt{
\mu} - \bigl(1+\epsilon^2\bigr)\sqrt{r}\bigr)_+^2,
\vspace*{3pt}
\cr
\displaystyle 1 - 2\alpha+\gamma> 2\bigl(\sqrt{\mu} - \bigl(1-
\epsilon^2\bigr)\sqrt{r}\bigr)_+^2 - \bigl(\sqrt{\mu} -
\epsilon^2\sqrt{r}\bigr)_+^2. } %
\end{equation}
In order to find appropriate $\mu\in(0, 1- \gamma)$ and $\epsilon
>0$ depending only on $\gamma, \alpha, r$ such that these two
inequalities hold simultaneously, we will discuss three cases separately.

\begin{longlist}
\item[\textit{Case} 1: $\frac{1+\gamma}{2} \leq\alpha< \frac{3
+\gamma}{4}$ and $\alpha- \frac{1+\gamma}{2} < r \leq\frac{1 -
\gamma}{4}$.]
In this case, let $\mu= 4r(1 - \epsilon^2)^2$. Then both inequalities
in \eqref{eqconstraints} hold when we let $\epsilon=0$. By the
continuity of the functions with respect to $\epsilon$ and the
properties of open sets, we can choose a sufficiently small positive
constant $C_0(\gamma, \alpha, r)$, such that when $\epsilon=
C_0(\gamma, \alpha, r) >0$, both inequalities in \eqref
{eqconstraints} hold strictly.\vspace*{1pt}

\item[\textit{Case} 2: $\frac{1+\gamma}{2} \leq\alpha< \frac{3
+\gamma}{4}$ and $r > \frac{1 - \gamma}{4}$.] In this case, let $\mu
= (1 - \gamma)(1 - \epsilon^2)^2$. Then both inequalities in \eqref
{eqconstraints} hold when we let $\epsilon=0$. Similarly, they also
hold when $\epsilon= C_0(\gamma, \alpha, r) > 0$.

\item[\textit{Case} 3: $\frac{3+\gamma}{4} \leq\alpha<1$ and $r>
(\sqrt{1 - \gamma} - \sqrt{ 1- \alpha} )^2$.] In this case,
let $\mu= (1 - \gamma)(1 - \epsilon^2)^2$. Then both inequalities in
\eqref{eqconstraints} hold when we let $\epsilon=0$. Similarly, they
also hold when $\epsilon= C_0(\gamma, \alpha, r) > 0$.

In summary, for fixed $\frac{1+\gamma}{2}\leq\alpha<1$ and $r >
\rho_\gamma^*(\alpha)$, we can choose $\epsilon>0$ and $\mu\in(0,
(1-\gamma)(1-\epsilon)^2]$ only depending on $\gamma, \alpha, r$
such that such that both inequalities in \eqref{eqconstraints} hold.
Notice that $t=\sqrt{\mu\log p}$ lies in the domain of $\mathrm{HC}(t)$, that
is, $ [1, \sqrt{\log\frac{N}{3}} ]$. Then when $p >
C(\alpha, \gamma, r)$, we have the inequality $\mathbb{P}
(\mathrm{HC}(\sqrt
{\mu\log p})\leq\log^2 N )\leq\frac{1}{\log^2 N}$, and
hence $\mathbb{P} (\mathrm{HC}^*\leq\log^2 N )\leq\frac{1}{\log
^2 N}$.
Since $C(\alpha, \gamma, r)$ is independent of the choice of $(\bolds
{\tau}, \bolds{\tilde{\beta}})$, we have
\[
\lim_{p \rightarrow\infty} \max_{(\bolds{\tau},
\bolds{\tilde{\beta}}) \in\Gamma(p, N, s, r)}
\mathbb{P}_{(\bolds
{\tau},
\bolds{\tilde{\beta}})} (H_0 \mbox{ is failed to reject} ) = 0.
\]
\end{longlist}

\subsubsection{Upper bound of $\mathrm{HC}^*$ under the null}
Under the null, for $m=1, \ldots, q$ we have $v_m=w_m$. By \eqref
{eqindependence}, for any $u=1, \ldots, L$, the variables $v_u,
v_{u+L},\break v_{u+2L}, \ldots, v_{u+(N-1)L}$ are i.i.d. standard complex
normal variables. This implies that $1_{\{|v_u|\geq t\}}, 1_{\{
|v_{u+L}|\geq t\}}, \ldots, 1_{\{|v_{u+(N - 1)L}|\geq t\}}$ are i.i.d.
Bernoulli random variables with parameter $\bar{\Psi}(t) = e^{-t^2}$.

By Lemma~\ref{thmVanzuijlen}, for any $u$, we have
\[
\mathbb{P} \biggl(\sup_{1 \leq t \leq\sqrt{\log\sfrac
{N}{3}}}\frac
{ \sum_{j=1}^{N} 1_{\{|v_{u + (j-1) L}|>t\}} - N \bar{\Psi
}(t)}{\log N\sqrt{ N \bar{\Psi}(t)(1- \bar{\Psi}(t))}} >
C_1 \biggr) < \frac{1}{\log^2 N},
\]
provided $N>N_0$. Therefore,
\begin{eqnarray*}
\mathrm{HC}^*&=&\sup_{1 \leq t \leq\sqrt{\log\sfrac{N}{3} }}\frac
{1}{\sqrt{L}}\sum
_{u=1}^{L} \frac{\sum_{j=1}^{N} 1_{\{|v_{u
+ (j-1)L}|>t\}} - N\bar{\Psi}(t)}{\sqrt{N \bar{\Psi}(t)(1- \bar
{\Psi}(t))}}
\\
&\leq & \frac{1}{\sqrt{L}}\sum_{u=1}^{L}
\sup_{1 \leq
t \leq\sqrt{\log\sfrac{N}{3} }} \frac{\sum_{j=1}^{N} 1_{\{|v_{u +
(j-1)L}|>t\}} - N\bar{\Psi}(t)}{\sqrt{N \bar{\Psi}(t)(1- \bar
{\Psi}(t))}}\leq\sqrt{L}C_1 \log
N
\end{eqnarray*}
with probability at least $1 - \frac{L}{\log^2 N}$. Since $L=\lfloor
\log N + 1 \rfloor$ and $N=p^{1 - \gamma}$, we have
\[
\lim_{p \rightarrow\infty}\mathbb{P}(H_0 \mbox{ is rejected})=
\lim_{p
\rightarrow\infty} \mathbb{P}\bigl(\mathrm{HC}^* > \log^2 N\bigr)=0.
\]




\begin{supplement}[id=suppA]
\stitle{$\!$Supplement to ``Global testing against sparse alternatives in time-frequency analysis''}
\slink[doi]{10.1214/15-AOS1412SUPP} 
\sdatatype{.pdf}
\sfilename{aos1412\_supp.pdf}
\sdescription{We give in \cite{CEL2015} the proofs to Lemmas \ref
{teocomplextail}, \ref{thmindcorr} and Theorem~\ref{lowerboundsparse}.}
\end{supplement}






\printaddresses
\end{document}